\documentclass[10pt]{article}
\usepackage{amssymb,amsmath,amsthm}
\usepackage{indentfirst}
\usepackage{exscale}
\usepackage{relsize}
\usepackage{subfigure}
\usepackage{graphicx}

\newcommand{\R}{\mathbb{R}}
\theoremstyle{definition}

\theoremstyle{remark}

\numberwithin{equation}{section}
\allowdisplaybreaks
 \UseRawInputEncoding
\begin{document}
\title{\Large\bf{ Existence of solutions for a poly-Laplacian system involving concave-convex nonlinearities on locally finite graphs}}
\date{}
\author {Ping Yang$^1$, \ Xingyong Zhang$^{1,2}$\footnote{Corresponding author, E-mail address: zhangxingyong1@163.com}\\
{\footnotesize $^1$Faculty of Science, Kunming University of Science and Technology,}\\
 {\footnotesize Kunming, Yunnan, 650500, P.R. China.}\\
{\footnotesize $^{2}$Research Center for Mathematics and Interdisciplinary Sciences, Kunming University of Science and Technology,}\\
 {\footnotesize Kunming, Yunnan, 650500, P.R. China.}\\
      }

 \date{}
 \maketitle

\begin{center}
\begin{minipage}{15cm}
\par
\small  {\bf Abstract:} We investigate the existence of two nontrivial solutions for a poly-Laplacian system involving concave-convex nonlinearities and parameters with Dirichlet boundary  condition on  locally finite graphs.
By using the mountain pass theorem and Ekeland's variational principle,
we obtain that system has at least one nontrivial solution of positive energy and one  nontrivial  solution of negative energy, respectively.
We also obtain an estimate about semi-trivial solutions.
Moreover, by using a result in \cite{Brown 09} which is based on the fibering  method and Nehari manifold, we obtain the existence of ground state solution to the single equation corresponding to  poly-Laplacian system.
Especially, we present some ranges of parameters in all of results.
\par
{\bf Keywords:} poly-Laplacian system (equation), concave-convex nonlinearities, mountain pass theorem, Ekeland's variational principle,
 locally finite graph.
\par
 {\bf 2020 Mathematics Subject Classification:} 35J60; 35J62; 49J35.
\end{minipage}
 \end{center}

\section{Introduction }
\setcounter{equation}{0}
\par
Existence and multiplicity of solutions for elliptic partial differential equations and systems involving concave-convex nonlinearities in Euclidean space have attracted some attentions.
In \cite{Ambrosetti 1994}, Ambrosetti and Brezis studied a class of second order Laplacian equations involving concave-convex nonlinearities with constant coefficient. By using the sub- and supersolutions, and variational arguments, they obtained some existence and multiplicity results of solutions.
In \cite{Brown 2007}, Brown and Wu also studied the second order Laplacian equation involving concave-convex nonlinearities with weight functions. By using the fibering  method and the Nehari manifold which was introduced by Pohozaev in \cite{POHOZAEV}, they obtained that the equation has at least two nontrivial solutions.
Moreover, in \cite{Brown 09}, Brown and Wu studied a potential operator equation. By using the methods similar to \cite{Brown 2007}, they obtained that the equation has at least two nontrivial solutions when the functionals related to  potential operators satisfy some appropriate conditions.
In \cite{Chen 2011}, Chen-Kuo-Wu studied a class of second order Kirchhoff equation involving concave-convex nonlinearities and parameters. By using the fibering method and the Nehari manifold, they obtained the existence of multiple positive solutions.
In \cite{Chen 2013}, Chen-Huang-Liu studied a class of nonhomogeneous $p$-Kirchhoff equation involving concave-convex nonlinearities with weight functions and a perturbation. By using the mountain pass theorem and Ekeland's variational principle, they obtained that the equation has at least two nontrivial solutions.
In \cite{Wu 08}, Wu investigated the multiplicity of nontrivial solutions for the following second order Laplacian elliptic system:
\begin{eqnarray}\label{p1}
 \begin{cases}
  -\Delta u=\lambda f(x) |u|^{\gamma-2}u+\frac{\alpha}{\alpha+\beta} h(x)|u|^{\alpha-2}u|v|^{\beta},\ \ & x \in \Omega,\\
  -\Delta v=\mu g(x) |v|^{\gamma-2}v+\frac{\beta}{\alpha+\beta} h(x)|u|^{\alpha}|v|^{\beta-2}v,\ \ & x \in \Omega,\\
  u(x)=v(x)=0,\ \ & x \in \partial\Omega,\\
 \end{cases}
\end{eqnarray}
where $\Omega \subseteq \R^{N}$ is a bounded domain and $1<\gamma<2<\alpha+\beta<\frac{2N}{N-2}$, and $f,g$ are sign-changing weight functions.
By using the fibering method and the Nehari manifold, the existence results of two nontrivial nonnegative solutions are obtained when the pair of the parameters $(\lambda, \mu)$ belongs to a certain subset of $\mathbb R$.
In \cite{E-S 2023},  Echarghaoui-Sersif investigated the existence of infinitely many solutions for a class of second order semilinear elliptic system involving critical Sobolev growth
and concave nonlinearities.
 By using the Fountain theorem, they obtained that the system has two distinct and infinite sets of radial solutions. One of the set exhibits positive energy, and the other set exhibits negative energy.
In \cite{B-M 04}, Bozhkova-Mitidieri considered a class of $(p,q)$-Laplacian elliptic system with Diriclet boundary condition.
By using the fibering method and the Nehari manifold, they obtained the existence and multiplicity of solutions of system. Moreover, by using the Pohozaev identity, they obtained the nonexistence result of solutions.
In \cite{Liu-Ou 17}, Liu and Ou investigated the following $(p,q)$-Laplacian elliptic system:
\begin{eqnarray*}
 \begin{cases}
  -\Delta_p u=\lambda \alpha a(x) |u|^{\alpha-2}u|v|^{\beta}+ \gamma b(x)|u|^{\gamma-2}u|v|^{\eta},\ \ & x \in \Omega,\\
  -\Delta_q v=\lambda \beta a(x) |u|^{\alpha}u|v|^{\beta-2}+ \eta b(x)|u|^{\gamma}u|v|^{\eta-2},\ \ & x \in \Omega,\\
  u(x)=v(x)=0,\ \ & x \in \partial\Omega,\\
 \end{cases}
\end{eqnarray*}
where $\Omega \subseteq \R^{N}$ is a bounded domain, $1<p,q<N$, $\alpha,\beta,\gamma,\eta>0$, $1<\alpha+\beta<\min\{p,q\}$ and $\max\{p,q\}<\gamma+\eta<\min\{p^*,q^*\}$.
By using the fibering method and the Nehari manifold, they obtained that  system has at least two nontrivial solutions.
We refer to \cite{A-H 06, Hsu 09, Shao 2018, Chen 2015} for more results about the elliptic systems involving concave-convex nonlinearities.
\par
Moreover, in recent years, the researches about  equations on graph have also attracted some attentions. For example, see \cite{Yamabe 2016, Shao 2021,Han 2019}.
In \cite{Yamabe 2016}, Grigor'yan-Lin-Yang considered the second order Laplace equation with the nonlinear term satisfying supequadratic condition and some additional assumptions on finite graphs and locally finite graphs. They obtained that the equation has a nontrivial solution by using the mountain pass theorem.
Furthermore, they also investigated the $p$-Laplacian equation and poly-Laplacian equation on finite graphs and locally finite graphs and obtained some similar results.
In \cite{Shao 2021}, Han and Shao studied the $p$-Laplacian equation with Dirichlet boundary condition on the locally finite graph. They proved that the equation admits a positive solution and a ground state solution by using the mountain pass theorem and the method of Nehari manifold.
In addition, in \cite{Han 2019}, Han and Shao studied a nonlinear biharmonic equation with a parameter and Dirichlet boundary condition on the locally finite graph. When the parameter is small enough, they obtained the existence and convergence of ground state solutions by using the mountain pass theorem and the method of Nehari manifold.
\par
Next, we recall some basic knowledge and notations of  locally finite graph, which is taken from \cite{Yamabe 2016, Shao 2021, Han 2019}.
Let $G=(V, E)$ be a graph, where $V$ and $E$ denote the vertex set and edge set, respectively. $xy$ denotes the edge connecting $x$ with $y$.
$(V,E)$ is defined as a locally finite graph if  there are only finite edges $xy \in E$ for any $x \in V$, and $(V,E)$ is defined as  finite graph if  both $V$ and $E$ are finite sets. $\omega_{xy}>0$ is defined as the weight on the edge $xy\in E$
 and it is assumed that $\omega_{xy}=\omega_{yx}$.
 Define the degree of a vertex $x\in V$ by $\mbox{deg}(x)=\sum_{y\thicksim x}\omega_{xy}$, where $y\thicksim x$ denotes that $y \in V$ satisfying the edge $xy\in E$.
 The distance of two vertices $x,y$ is denoted by $d(x,y)$, and is defined as the minimal number of edges that connect $x$ with $y$.
 Let $\Omega\subset V$. If  $d(x,y)$ is uniformly bounded for any $x,y\in \Omega$, then $\Omega$ is called as a bounded domain in $V$. Define the boundary of $\Omega$ by
$$
\partial \Omega=\{y \in V,\; y \notin \Omega\big| \; \exists \; x \in \Omega \mbox{ such that } xy \in E\}.
$$
Let $\mu:V\to \R^+$ be a finite measure satisfying $\mu(x)\geq \mu_0>0$. For any function $u:V\to\R$, one denotes
\begin{eqnarray}
\label{eq12}
\int\limits_V u(x) d\mu=\sum\limits_{x\in V}u(x)\mu(x).
\end{eqnarray}
Let $C(V)=\{u|u:V\to\R\}$. Define the Laplacian operator $\Delta: C(V)\to C(V)$ by
\begin{eqnarray}
\label{eq8}
\Delta u(x)=\frac{1}{\mu(x)}\sum_{y\thicksim x}\omega_{xy}(u(y)-u(x))
\end{eqnarray}
and the associate gradient form is defined by
\begin{eqnarray}
\label{eq9}
\Gamma(u_1,u_2)(x)=\frac{1}{2\mu(x)}\sum\limits_{y\thicksim x}w_{xy}(u_1(y)-u_1(x))(u_2(y)-u_2(x)).
\end{eqnarray}
Denote $\Gamma(u):=\Gamma(u,u)$. The length of the gradient is defined by
\begin{eqnarray}
\label{eq10}
|\nabla u|(x)=\sqrt{\Gamma(u)(x)}=\left(\frac{1}{2\mu(x)}\sum\limits_{y\thicksim x}w_{xy}(u(y)-u(x))^2\right)^{\frac{1}{2}},
\end{eqnarray}
and the length of $m_{th}$ order gradient is defined by
\begin{equation}
\label{eq13}
|\nabla^m u|=
\begin{cases}
|\nabla \Delta^{\frac{m-1}{2}}u|, \ \ \mbox{if $m$ is odd},\\
|\Delta ^{\frac{m}{2}} u|,\ \ \mbox{if $m$ is even}.
\end{cases}
\end{equation}
For any given $s> 1$,  the $s$-Laplacian operator $ \Delta_{s}:C(V)\rightarrow C(V) $ is defined as follows:
\begin{equation}
\label{eq14}
  \Delta_{s}u(x) = \frac{1}{2\mu(x)} \sum_{y\sim x}  \left(|\nabla u|^{s-2} (y)+ |\nabla u|^{s-2} (x) \right) \omega_{xy} (u(y)-u(x)).
  \end{equation}
Let $C_c(\Omega):=\{u:V \to \R \big| supp\mbox{ }u\subset \Omega \mbox{ and }\forall x \in V \backslash\Omega,u(x)=0\}$. For any function $\phi \in C_c(\Omega)$, there holds
\begin{eqnarray}
\label{eq11}
\int_{\Omega} \Delta_{s}u\phi d\mu=-\int_{\Omega\cup \partial\Omega} |\nabla u|^{s-2}\Gamma(u,\phi)d\mu.
\end{eqnarray}
 For any $1\leq r<+\infty$, let $L^r(\Omega)$ be the completion of $C_c(\Omega)$ under the norm
$$
\|u\|_{L^r(\Omega)}=
\left(\int_{\Omega}|u(x)|^r d \mu\right)^{\frac{1}{r}}.
$$
Let $W_0^{m,s}(\Omega)$ be the completion of $C_c(\Omega)$ under the norm
$$
\|u\|_{W_0^{m,s}(\Omega)}=\left(\int_{\Omega\cup\partial\Omega}|\nabla^m u(x)|^sd\mu\right)^\frac{1}{s},
$$
where $m$ is a positive integer and  $s>1$. For any $u\in W_0^{m,s}(\Omega)$, we also define the following norm:
$$
\|u\|_{\infty}=\max_{x\in \Omega}|u(x)|.
$$
$W_0^{m,s}(\Omega)$ is of finite dimension. See \cite{Yamabe 2016, Shao 2021} for more details.

\par
In this paper, our work are mainly inspired by \cite{Brown 09, Chen 2013, Wu 08, Liu-Ou 17, Yamabe 2016}. We shall apply the mountain pass theorem and Ekeland's variational principle like \cite{Chen 2013} to investigate the multiplicity of solutions for a class of poly-Laplacian systems on graph, which can be seen as a discrete version of (\ref{p1}) on graph in some sense,  and we also apply an abstract result in \cite{Brown 09} which is essentially obtained by the fibering  method and the Nehari manifold like \cite{Brown 2007, Brown 09, Wu 08,Liu-Ou 17} to investigate the multiplicity of solutions for a class of poly-Laplacian equations on graph. To be specific,
we  investigate the existence of nontrivial solutions for the following poly-Laplacian system on a locally finite graph $G=(V,E)$:
\begin{eqnarray}
\label{eq1}
 \begin{cases}
  \pounds_{m_1,p}u=\lambda_1 h_1(x)|u|^{\gamma_1-2}u+\frac{\alpha}{\alpha+\beta}c(x)|u|^{\alpha-2}u|v|^{\beta},&\;\;\;\;\hfill x\in \Omega,\\
   \pounds_{m_2,q}v=\lambda_2 h_2(x)|v|^{\gamma_2-2}v+\frac{\beta}{\alpha+\beta}c(x)|u|^{\alpha}|v|^{\beta-2}v,&\;\;\;\;\hfill x\in \Omega,\\
   u=v=0,&\;\;\;\;\hfill x\in \partial\Omega,
   \end{cases}
\end{eqnarray}
where $\Omega\cup \partial\Omega\subset V$ is a bounded domain, $m_i,  i=1,2$ are positive integers, $p,q,\gamma_1,\gamma_2>1$, $\lambda_1,\lambda_2,\alpha,\beta>0$, $\max\{\gamma_1,\gamma_2\} < \min\{p,q\} \le \max\{p,q\} < \alpha+\beta$, $h_1(x),h_2(x), c(x):\Omega \to \R^+$, and $\pounds_{m_i,s}$ ($i=1,2, s=p,q)$ are defined as follows: for any function $\phi: \Omega \cup \partial\Omega \to \R$,
\begin{eqnarray}
\label{eq15}
\int_{ \Omega} (\pounds_{m_i,s}u) \phi d\mu=
\begin{cases}
\int_{\Omega \cup \partial\Omega} |\nabla^{m_i} u|^{s-2}\Gamma(\Delta^{\frac{m_i-1}{2}}u,\Delta^{\frac{m_i-1}{2}}\phi), \ \ \mbox{if $m_i$ is odd},\\
\int_{\Omega \cup \partial\Omega}|\nabla^{m_i} u|^{s-2}\Delta ^{\frac{m_i}{2}} u\Delta ^{\frac{m_i}{2}} \phi,\ \ \mbox{if $m_i$ is even}.
\end{cases}
\end{eqnarray}
When $m=1$, $\pounds_{m,p}u=-\Delta_p u$, and when $p=2$, $\pounds_{m,p}u=(-\Delta^m)u$ which is called as poly-Laplacian operator of $u$.  More details can be seen in  \cite{Yamabe 2016} for the definition of $\pounds_{m,p}$. It is easy to see that the system (\ref{eq1}) with $m=1$, $p=2$ and $\gamma_1=\gamma_2=\gamma$ is a generalization of (\ref{p1}) from the Euclidean  setting to locally finite graph.

\par
In this paper, we call that $(u,v)$ is a semi-trivial solution of system (\ref{eq1}) if $(u,v)$ is a solution of system (\ref{eq1}) with either $(u,v)=(u,0)$ or $(u,v)=(0,v)$, and we call that $(u,v)$ is a nontrivial solution of system (\ref{eq1}) if $(u,v)$ is a solution of system (\ref{eq1}) and $(u,v)\not=(0,0)$.

\par

Denote
\begin{eqnarray*}
M_{(\lambda_1,\lambda_2)} & = & 2^{1-\max\{p,q\}}\min\left\{ \frac{1-\lambda_1 C_{m_1,p}^p(\Omega)}{p},\frac{1- \lambda_2C_{m_2,q}^q(\Omega)}{q}\right\},\\
M_2 & = & \frac{C_0}{(\alpha+\beta)^2} \left(\alpha C_{m_1,p}^{\alpha+\beta}(\Omega) + \beta C_{m_2,q}^{\alpha+\beta}(\Omega) \right),
\end{eqnarray*}
where $C_0=\max_{x\in \Omega}c(x)$ and $C_{m_1,p}(\Omega)$ and $C_{m_2,q}(\Omega)$ are embedding constants given in Lemma 2.1 below. Especially, we present a concrete value of $C_{m_1,p}(\Omega)$ and $C_{m_2,q}(\Omega)$ if $m_1=m_2=1$, $p,q\geq2$ and for any $x \in \Omega$, there exists at least one $y \in \partial \Omega$ such that $y \thicksim x$ (see Lemma 2.2 below).

\par

Next, we state our main results. Assume that $\lambda_1$ and $\lambda_2$ satisfy the following inequalities:
{\footnotesize\begin{eqnarray}
\begin{cases}
\label{eqq0}
 &  0 < \lambda_1 < C_{m_1,p}^{-p}(\Omega),\\
 &  0 < \lambda_2 < C_{m_2,q}^{-q}(\Omega),\\
 &  M_{(\lambda_1,\lambda_2)} \le \frac{\alpha+\beta}{\max\{p,q\}}M_2,\\
 & \frac{\lambda_1(p-\gamma_1)} {p\gamma_1}\|h_1\|_{L^{\frac{p}{p-\gamma_1}}(\Omega)}^{\frac{p}{p-\gamma_1}}
    +  \frac{\lambda_2(q-\gamma_2)}{q\gamma_2}\|h_2\|_{L^{\frac{q}{q-\gamma_2}}(\Omega)}^{\frac{q}{q-\gamma_2}}
   <  \frac{\alpha+\beta - \max\{p,q\}}{\alpha+\beta}
        M_{(\lambda_1,\lambda_2)}^{\frac{\alpha+\beta}{\alpha+\beta-\max\{p,q\}}}
        \left(\frac{\max\{p,q\}}{(\alpha+\beta)M_2}\right)^{\frac{\max\{p,q\}}{\alpha+\beta-\max\{p,q\}}}.
\end{cases}
\end{eqnarray}
}
\vskip2mm
\noindent
{\bf Theorem 1.1.} {\it Let $G=(V,E)$ be a locally finite graph, $\Omega\not=\emptyset$ and $\partial \Omega\not=\emptyset$. If $(\lambda_1,\lambda_2)\in (0,+\infty)\times(0, +\infty)$ satisfies (\ref{eqq0}), system (\ref{eq1}) has at least one nontrivial solution of  positive energy and one nontrivial  solution of negative energy.
 }
\vskip2mm
\noindent
{\bf Remark 1.1.} {\it There exist $\lambda_1$ and $\lambda_2$ satisfying (\ref{eqq0}). For example,
let $m_1=2$, $m_2=3$, $\gamma_1=2$, $\gamma_2=3$, $p=4$, $q=5$, $\alpha=2$, $\beta=4$, and
$$
C_0  =  \frac{1}{C_{2,4}^6(\Omega)+2C_{3,5}^6(\Omega)},\ \ \|h_1\|_{L^2(\Omega)}^2  =  \frac{5^{6}}{9\cdot2^{31}}C_{2,4}^4(\Omega),\ \ \|h_2\|_{L^{\frac{5}{2}}(\Omega)}^{\frac{5}{2}}  =  \frac{5^{6}}{2^{32}}C_{3,5}^5(\Omega).
$$
When $\lambda_1=\frac{1}{5}C_{2,4}^{-4}(\Omega)$ and $\lambda_2=\frac{1}{6}C_{3,5}^{-5}(\Omega)$,
we can obtain that
$$
M_{(\lambda_1,\lambda_2)}=\frac{1}{3\cdot2^5}=\frac{1}{96},\ \  M_2=\frac{1}{18}.
$$
Evidently,
$$
M_{(\lambda_1,\lambda_2)} < \frac{6}{5}M_2.
$$
Moreover,
$$
  \frac{1}{20}C_{2,4}^{-4}(\Omega)\|h_1\|_{L^2(\Omega)}^2+\frac{1}{45}C_{3,5}^{-5}(\Omega)\|h_2\|_{L^{\frac{5}{2}}(\Omega)}^{\frac{5}{2}}
= \frac{5^5}{9\cdot2^{33}}+\frac{5^5}{9\cdot2^{32}}
< \frac{5^5}{9\cdot{2}^{31}}.
$$
Hence, (\ref{eqq0}) holds for $\lambda_1=\frac{1}{5}C_{2,4}^{-4}(\Omega)$ and $\lambda_2=\frac{1}{6}C_{3,5}^{-5}(\Omega)$.
}
\vskip2mm
\noindent
{\bf Theorem 1.2.} {\it Let $G=(V,E)$ be a locally finite graph, $\Omega\not=\emptyset$ and $\partial \Omega\not=\emptyset$. For each $\lambda_1>0$, if $(u,v)$ is a semi-trivial solution of  system (\ref{eq1}) and $(u,v)=(u,0)$, then
 $$
 \|u\|_{W^{m_1,p}_0(\Omega)} \le \left(\lambda_1 H_1C^{\gamma_1}_{m_1,p}(\Omega)\right)^{\frac{1}{p-\gamma_1}},
 $$
where $H_1=\max_{x\in \Omega} h_1(x)$.
 Similarly, for each $\lambda_2>0$,  if $(u,v)$ is a semi-trivial solution of  system (\ref{eq1}) and $(u,v)=(0,v)$, then
 $$
 \|v\|_{W^{m_2,q}_0(\Omega)} \le \left(\lambda_2 H_2C^{\gamma_2}_{m_2,q}(\Omega)\right)^{\frac{1}{q-\gamma_2}},
 $$
  where   $H_2=\max_{x\in \Omega} h_2(x)$.
 }

\par
Moreover, we also investigate the existence of ground state solution for the following poly-Laplacian equation on  a locally finite graph $G=(V,E)$ by applying  Theorem 3.3 in \cite{Brown 09}:
\begin{eqnarray}
\label{a1}
\begin{cases}
 \pounds_{m,p}u=\lambda h(x)|u|^{\gamma-2}u+c(x)|u|^{\alpha-2}u,&\;\;\;\;\hfill x\in \Omega,\\
 u(x)=0,&\;\;\;\;\hfill x\in \partial\Omega,
\end{cases}
\end{eqnarray}
where $\Omega\cup \partial\Omega\subset V$ is a bounded domain, $m$ is a positive integer, $p, \gamma>1$, $\lambda>0$, $\gamma <p <\alpha$ and $h(x),c(x):\Omega \to \R^+$.
Denote
\begin{eqnarray*}
\lambda_0
= \frac{p-\gamma}{H_0} C_{m,p}^{-\gamma}(\Omega) \left(\left(C_0 C_{m,p}^{\alpha}(\Omega)\right)^{p-\alpha}(\alpha-p)^{\alpha-p}(\alpha-\gamma)^{\gamma-\alpha}\right)^{\frac{1}{p-\gamma}},
\end{eqnarray*}
 \begin{eqnarray} \label{aa2}
  \lambda_{\star}
=  \frac{\gamma(\alpha-p)}{p\alpha H_0 C_{m,p}^\gamma(\Omega)} \left(C_0C^{\alpha}_{m,p}(\Omega)\right)^{\frac{p- \gamma}{p-\alpha}},\ \
\lambda_{\star\star}
=  \min\{\lambda_0, \lambda_{\star}\},
\end{eqnarray}
where $H_0=\max_{x \in \Omega} h(x)$ and $C_0=\max_{x \in \Omega} c(x)$. We obtain the following result.
\vskip2mm
\noindent
{\bf Theorem 1.3.} {\it Let $G=(V,E)$ be a locally finite graph, $\Omega\not=\emptyset$ and $\partial \Omega\not=\emptyset$. If $\lambda \in (0, \lambda_0)$, then equation (\ref{a1}) has at least one nontrivial solution of positive energy and one nontrivial solution of negative energy. Furthermore, if $\lambda \in (0, \lambda_{\star\star})$, the solution of negative energy is the ground state solution of (\ref{a1}).}

\vskip2mm
\noindent
{\bf Remark 1.2.}  Similar to Theorem 1.1, we can also obtain  one nontrivial solution of positive energy and one nontrivial solution of negative energy for  (\ref{a1}) by using the mountain pass theorem and Ekeland's variational principle. We do not know whether these two solutions are different from those two solutions in Theorem 1.3 which is obtained essentially by the fibering method and Nehari manifold.

\section{Preliminaries}
\setcounter{equation}{0}
\par
Define the space $W=W_0^{m_1,p}(\Omega)\times W_0^{m_2,q}(\Omega)$ with the norm
$$
\|(u,v)\|_W=\|u\|_{W_0^{m_1,p}(\Omega)}+\|v\|_{W_0^{m_2,q}(\Omega)}.
$$
Then $W$ is a Banach space and of finite dimension.
 Consider the energy functional $\varphi:\; W \to \R$ defined as
\begin{eqnarray}
\label{eq21}
          \varphi(u,v)
 &  =  &  \frac{1}{p}\int_{\Omega \cup \partial\Omega} |\nabla^{m_1} u|^pd\mu-\frac{\lambda_1}{\gamma_1}\int_{\Omega}h_1(x)|u|^{\gamma_1}d\mu\nonumber\\
 &     &  + \frac{1}{q}\int_{\Omega \cup \partial\Omega} |\nabla^{m_2} v|^qd\mu-\frac{\lambda_2}{\gamma_2}\int_{\Omega}h_2(x)|v|^{\gamma_2}d\mu-\frac{1}{\alpha+\beta}\int_{\Omega} c(x)|u|^{\alpha}|v|^{\beta}d\mu.
\end{eqnarray}
Then $\varphi(u,v) \in C^1(W,\R)$ and
\begin{eqnarray}
\label{eq211}
           \langle \varphi'(u,v),(\phi,\psi)\rangle
 &  =  &   \int_{\Omega \cup \partial\Omega} (\pounds_{m_1,p}u) \phi d\mu-\lambda_1\int_{\Omega}h_1(x)|u|^{\gamma_1-2}u\phi d\mu\nonumber\\
 &     & + \int_{\Omega \cup \partial\Omega} (\pounds_{m_2,q}v) \psi d\mu-\lambda_2\int_{\Omega}h_2(x)|v|^{\gamma_2-2}\psi d\mu\nonumber\\
 &     & - \frac{\alpha}{\alpha+\beta}\int_{\Omega} c(x)|u|^{\alpha-2}u|v|^{\beta}\phi d\mu -\frac{\beta}{\alpha+\beta}\int_{\Omega} c(x)|u|^{\alpha}|v|^{\beta-2}v \psi d\mu.
\end{eqnarray}

\vskip2mm
\noindent
{\bf Definition 2.1. }{\it  $(u,v)\in W$ is called as a weak solution of system (\ref{eq1}) if
\begin{eqnarray}
 \label{eq22}
 &   &       \int_{\Omega \cup \partial\Omega}(\pounds_{m_1,p}u) \phi d\mu
         =   \lambda_1\int_{\Omega}h_1(x)|u|^{\gamma_1-2}u\phi d\mu +\frac{\alpha}{\alpha+\beta}\int_{\Omega} c(x)|u|^{\alpha-2}u|v|^{\beta}\phi d\mu,\\
\label{eq23}
&    &       \int_{\Omega \cup \partial\Omega}(\pounds_{m_2,q}v) \psi d\mu
         =   \lambda_2\int_{\Omega}h_2(x)|v|^{\gamma_2-2}v\psi d\mu+\frac{\beta}{\alpha+\beta}\int_{\Omega} c(x)|u|^{\alpha}|v|^{\beta-2}v \psi d\mu,
\end{eqnarray}
for all $(\phi,\psi)\in W$.
}
\par
Evidently, $(u,v)\in W$ is  a weak solution of system (\ref{eq1}) if and only if $(u,v)$ is a critical point of $\varphi$.
\vskip2mm
\noindent
{\bf Proposition 2.2. }{\it If  $(u,v)\in W$ is a weak solution of system (\ref{eq1}), then $(u,v)\in W$ is also a point-wise solution of (\ref{eq1}). }
\vskip2mm
\noindent
{\bf Proof.} For any
fixed $y\in V$, we take a test function $\phi:V\to \R$ in (\ref{eq22}) with
\begin{eqnarray*}
\phi(x)=\begin{cases}
1,& x=y,\\
0,& x\not=y,
\end{cases}
\end{eqnarray*}
and a test function $\psi:V\to \R$ in (\ref{eq23}) with
\begin{eqnarray*}
\psi(x)=\begin{cases}
1,& x=y,\\
0,& x\not=y.
\end{cases}
\end{eqnarray*}
Thus, we have
\begin{eqnarray*}
&   &    \pounds_{m_1,p}u(y)  =  \lambda_1 h_1(y)|u(y)|^{\gamma_1-2}u(y) + \frac{\alpha}{\alpha+\beta}c(y)|u(y)|^{\alpha-2}u(y)|v(y)|^{\beta},\\
 &  &     \pounds_{m_2,q}v(y)  =  \lambda_2 h_2(y)|v(y)|^{\gamma_2-2}v(y) + \frac{\beta}{\alpha+\beta} c(y)|u(y)|^{\alpha}|v(y)|^{\beta-2}v(y).
\end{eqnarray*}
By the arbitrary of $y$, we complete the proof.\qed
\vskip2mm
\noindent

\vskip2mm
 \noindent
{\bf Lemma 2.1.} (Embedding theorem \cite{Yamabe 2016}). {\it Let $G=(V, E)$ be a locally finite graph, $\Omega\cup \partial\Omega\subset V$ be a bounded domain with $\Omega\not=\emptyset$. For any $m \ge 1,s>1$, $W_0^{m,s}(\Omega)$ is embedded in
$L^r(\Omega)$ for all $1 \le r \le+\infty$. In particular, there exists a constant $C_{m,s}(\Omega)>0$ depending only on $m,s$ and $\Omega$ such that
\begin{eqnarray}
\label{eq24}
\left(\int_{\Omega}|u(x)|^r d \mu\right)^{\frac{1}{r}} \le C_{m,s}(\Omega)\left(\int_{\Omega\cup\partial\Omega}|\nabla^m u(x)|^sd\mu\right)^\frac{1}{s},
\end{eqnarray}
where
\begin{eqnarray}
\label{eqc0}
C_{m,s}(\Omega) = \frac{C}{\mu_{min}}(1+|\Omega|)
\ \ \mbox{with}\ \  C \ \ \mbox{satisfying}\ \
\|u\|_{L^r(\Omega)} \le C \|u\|_{W_0^{m,r}(\Omega)},
\end{eqnarray}
and  $|\Omega|=\sum_{x \in \Omega} \mu(x)$.
Moreover, $W_0^{m,s}(\Omega)$ is pre-compact, namely, if $\{u_k\}$ is bounded in $W_0^{m,s}(\Omega)$, then up to a subsequence, there exists some $u \in W_0^{m,s}(\Omega)$ such that $u_k \rightarrow u$ in $W_0^{m,s}(\Omega)$ as $k \rightarrow \infty$.}

\vskip2mm
\par
Furthermore, if for any $x \in \Omega$, there exists at least one $y \in \partial \Omega$ such that $y \thicksim x$,
we can present a specific value of $C_{m,s}(\Omega)$ with $m=1$ and $s\geq 2$. In details, we obtain the following lemma.
\vskip2mm
 \noindent
{\bf Lemma 2.2.} {\it Let $G=(V, E)$ be a locally finite graph. $\Omega\cup \partial\Omega\subset V$ is a bounded domain with $\Omega\not=\emptyset$ and $\partial\Omega\not=\emptyset$, and   for any $x \in \Omega$, there exists at least one $y \in \partial \Omega$ such that $y \thicksim x$. Then for any $s\geq2$ and $1 \le r <+\infty$,
\begin{eqnarray}
\label{224}
\left(\int_{\Omega}|u(x)|^r d \mu\right)^{\frac{1}{r}} \le C_{1,s}(\Omega)\left(\int_{\Omega\cup\partial\Omega}|\nabla u(x)|^sd\mu\right)^\frac{1}{s},
\end{eqnarray}
where
\begin{eqnarray*}
\label{eq28}
 C_{1,s}(\Omega) = (1+|\Omega|) \hat{\mu}_{\min}^{-\frac{1}{s}} \left(\frac{2\mu_{\max}}{w_{\min}}\right)^{\frac{1}{2}},
\end{eqnarray*}
$\hat{\mu}_{\min}=\min_{x\in \Omega}\mu(x)$, $\mu_{\max}=\max_{x\in \Omega\cup\partial \Omega}\mu(x)$ and $w_{\min}=\min_{x\in \Omega\cup\partial \Omega}w_{xy}$.
}
\vskip2mm
\noindent
{\bf Proof.} There holds
\begin{eqnarray*}
\label{eq25}
         \|u\|_{W_0^{1,s}(\Omega)}^s
&  =  &  \int_{\Omega\cup \partial \Omega} |\nabla u(x)|^sd\mu \nonumber\\
&  =  &  \sum_{x \in \Omega\cup \partial \Omega} \left(\frac{1}{2\mu(x)}\sum\limits_{y\thicksim x}w_{xy}(u(y)-u(x))^2\right)^{\frac{s}{2}} \mu(x)\nonumber\\
& \ge &  \left(\frac{w_{\min}}{2\mu_{\max}}\right)^{\frac{s}{2}}  \sum_{x \in \Omega\cup \partial \Omega} \sum\limits_{y\thicksim x}|u(y)-u(x)|^s \mu(x) \nonumber\\
&  =  &  \left(\frac{w_{\min}}{2\mu_{\max}}\right)^{\frac{s}{2}}  \left(\sum_{x \in \Omega}\sum\limits_{y\thicksim x}|u(y)-u(x)|^s  +  \sum_{x \in \partial \Omega} \sum\limits_{y\thicksim x}|u(y)-u(x)|^s\right) \mu(x) \nonumber\\
& \ge &  \left(\frac{w_{\min}}{2\mu_{\max}}\right)^{\frac{s}{2}} \left(\sum_{x \in \Omega}\sum\limits_{y\thicksim x,y\in \Omega}|u(y)-u(x)|^s  +  \sum_{x \in \Omega} \sum\limits_{y\thicksim x,y\in \partial \Omega}|u(x)|^s\right) \mu(x) \nonumber\\
& \ge &  \left(\frac{w_{\min}}{2\mu_{\max}}\right)^{\frac{s}{2}}  \sum_{x \in \Omega} \sum\limits_{y\thicksim x,y\in \partial \Omega}|u(x)|^s \mu(x) \nonumber\\
& \ge &  \left(\frac{w_{\min}}{2\mu_{\max}}\right)^{\frac{s}{2}}   \sum_{x \in \Omega} |u(x)|^s \mu(x).
\end{eqnarray*}
Hence,
\begin{eqnarray}
\label{eq25}
         \|u\|_{L^s(\Omega)} \le \left(\frac{2\mu_{\max}}{w_{\min}}\right)^{\frac{1}{2}} \|u\|_{W_0^{1,s}(\Omega)}.
\end{eqnarray}
Moreover, by (\ref{eq25}), we have
\begin{eqnarray*}
\label{eq26}
    \|u\|_{\infty}
\le \hat{\mu}_{\min}^{-\frac{1}{s}} \|u\|_{L^s(\Omega)}
\le \hat{\mu}_{\min}^{-\frac{1}{s}} \left(\frac{2\mu_{\max}}{w_{\min}}\right)^{\frac{1}{2}} \|u\|_{W_0^{1,s}(\Omega)},
\end{eqnarray*}
and for  any $1 \le r <+\infty$ , we have
\begin{eqnarray*}
\label{eq27}
         \|u\|_{L^r(\Omega)}
&  =  &  \left(\sum_{x \in \Omega} |u(x)|^r \mu(x) \right)^{\frac{1}{r}} \nonumber\\
& \le &  |\Omega|^{\frac{1}{r}} \|u\|_{\infty} \nonumber\\
& \le &  |\Omega|^{\frac{1}{r}} \hat{\mu}_{\min}^{-\frac{1}{s}} \left(\frac{2\mu_{\max}}{w_{\min}}\right)^{\frac{1}{2}} \|u\|_{W_0^{1,s}(\Omega)}\nonumber\\
& \le &  (1+|\Omega|) \hat{\mu}_{\min}^{-\frac{1}{s}} \left(\frac{2\mu_{\max}}{w_{\min}}\right)^{\frac{1}{2}} \|u\|_{W_0^{1,s}(\Omega)}.
\end{eqnarray*}
The proof is complete.
\qed

\vskip2mm
Let $X$ be a real Banach space. Suppose that $\varphi \in C^{1}(X, \mathbb{R})$. If any sequence $\{u_n\}$ which satisfies that $\varphi(u_n)$ is bounded for all $n\in \mathbb N$ and  $\varphi'(u_n)\to 0$ as $n\to \infty$  has a convergent subsequence, then we call that $\varphi$ satisfies the Palais-Smale condition.

\vskip2mm
\noindent
{\bf Lemma 2.3.} (Mountain pass theorem \cite{Rabinowitz 1986}) {\it Let $X$ be a real Banach space.  $\varphi \in C^{1}(X,\R)$, $\varphi(0)=0$ and $\varphi$ satisfies Palais-Smale condition. Suppose that $\varphi $ satisfies the following conditions:\\
(i) there exist two positive constants $\rho$ and $\alpha$ such that $ \varphi|_{\partial B_{\rho}(0)}\ge \alpha $, where $B_\rho=\{w\in X:\|w\|_X<\rho\}$;\\
(ii) there exists $ w\in X\backslash \bar B_{\rho} (0)$ such that $ \varphi(w)\le 0 $.\\
Then $\varphi$ has a critical value $c_*\ge \alpha$ with
 $$
 c_*:=\inf_{\gamma\in\Gamma}\max_{t\in[0,1]}\varphi(\gamma(t)),
 $$
where }
 $$
 \Gamma:=\{\gamma\in C([0,1],X]):\gamma(0)=0,\gamma(1)=w\}.
 $$

\vskip2mm
\noindent
{\bf Lemma 2.4.} (Ekeland's variational principle \cite{Mawhin J 1989}) {\it Let $X$ be a  complete metric space with metric $d$, and $\varphi:X\rightarrow \R$ be a lower semicontinuous function, bounded from below and not identical to $+\infty$. Let $\varepsilon >0$ be given and $U \in X$ such that
$$
\varphi(U)\le \inf_{X} \varphi+\varepsilon.
$$
Then there exists $V \in X$ such that
$$
\varphi(V)\le \varphi(U),\;d(U,V)\le1,
$$
and for each $O\in X$, one has}
$$
\varphi(V)\le \varphi(O)+\varepsilon d(V,O).
$$

\section{Proofs for Theorem 1.1 and Theorem 1.2}
\setcounter{equation}{0}
\vskip2mm
\noindent
{\bf Lemma 3.1.} {\it For each $(\lambda_1,\lambda_2)$ satisfying (\ref{eqq0}),
there exists a positive constant $\rho_{(\lambda_1,\lambda_2)}$ such that $\varphi(u,v)>0$ whenever $||(u,v)||_W=\rho_{(\lambda_1,\lambda_2)}$.}\\
{\bf Proof.} Note that
\begin{eqnarray}
\label{eq311}
c(x) \le \max_{x\in \Omega} c(x):=C_0,\ \ h_i^{\star}:=\min_{x\in \Omega} h_i(x) \le h_i(x) \le \max_{x\in\Omega}h_i(x):=H_i,\ \ i=1,2, \mbox{ for all } x\in \Omega.
\end{eqnarray}
Then, by Young's inequality and Lemma 2.1, we have
\begin{eqnarray}
\label{eq31}
&       &   \int_{\Omega}h_1(x)|u|^{\gamma_1}d\mu\nonumber\\
&  \le  &   \frac{p-\gamma_1}{p}\int_{\Omega}h_1(x)^{\frac{p}{p-\gamma_1}}d\mu+\frac{\gamma_1}{p}\int_{\Omega}|u|^pd\mu\nonumber\\
&  \le  &   \frac{p-\gamma_1}{p}\|h_1\|_{L^{\frac{p}{p-\gamma_1}}(\Omega)}^{\frac{p}{p-\gamma_1}}+\frac{\gamma_1}{p}C_{m_1,p}^p(\Omega) \|u\|_{W_0^{m_1,p}(\Omega)}^p,
\end{eqnarray}
for all $u\in W_0^{m_1,p}(\Omega)$. Similarly, we have
\begin{eqnarray}
\label{eq32}
       \int_{\Omega}h_2(x)|v|^{\gamma_2}d\mu
  \le   \frac{q-\gamma_2}{q}\|h_2\|_{L^{\frac{q}{q-\gamma_2}}(\Omega)}^{\frac{q}{q-\gamma_2}}+\frac{\gamma_2}{q}C_{m_2,q}^q(\Omega) \|v\|_{W_0^{m_2,q}(\Omega)}^q,
\end{eqnarray}
for all $v\in W_0^{m_2,q}(\Omega)$. Moreover,  it follows from  Young's inequality and Lemma 2.1 that
\begin{eqnarray}
\label{eq33}
&       &   \int_{\Omega}c(x)|u|^{\alpha}|v|^{\beta}d\mu\nonumber\\
&  \le  &   C_0 \int_{\Omega}|u|^{\alpha}|v|^{\beta}d\mu\nonumber\\
&  \le  &   C_0 \left(\frac{\alpha}{\alpha+\beta}\int_{\Omega}|u|^{\alpha+\beta} d\mu +\frac{\beta}{\alpha+\beta} \int_{\Omega}|v|^{\alpha+\beta}d\mu\right)\nonumber\\
&  \le  &   C_0\left(\frac{\alpha C_{m_1,p}^{\alpha+\beta}(\Omega)}{\alpha+\beta} \|u\|_{W_0^{m_1,p}(\Omega)}^{\alpha+\beta} + \frac{\beta C_{m_2,q}^{\alpha+\beta}(\Omega)}{\alpha+\beta}  \|v\|_{W_0^{m_2,q}(\Omega)}^{\alpha+\beta}\right),
\end{eqnarray}
for all $(u,v)\in W$.
Thus, (\ref{eq21}) and (\ref{eq31})-(\ref{eq33}) imply that when $(\lambda_1,\lambda_2) \in (0,C_{m_1,p}^{-p}(\Omega))\times (0,C_{m_2,q}^{-q}(\Omega))$, for any $(u,v) \in W$ with $\|(u,v)\|_W \le 1$, we have
\begin{eqnarray}
\label{eq34}
            \varphi(u,v)
&  =   &    \frac{1}{p} \|u\|_{W_0^{m_1,p}(\Omega)}^p - \frac{\lambda_1}{\gamma_1}\int_{\Omega}h_1(x)|u|^{\gamma_1}d\mu\nonumber\\
&      & +  \frac{1}{q} \|v\|_{W_0^{m_2,q}(\Omega)}^q - \frac{\lambda_2}{\gamma_2}\int_{\Omega}h_2(x)|v|^{\gamma_2}d\mu - \frac{1}{\alpha+\beta}\int_{\Omega} c(x)|u|^{\alpha}|v|^{\beta}d\mu\nonumber\\
& \ge  &    \frac{1}{p} \left(1-\lambda_1 C_{m_1,p}^p(\Omega)\right)\|u\|_{W_0^{m_1,p}(\Omega)}^p
         - \frac{\lambda_1(p-\gamma_1)} {p\gamma_1}\|h_1\|_{L^{\frac{p}{p-\gamma_1}}(\Omega)}^{\frac{p}{p-\gamma_1}}\nonumber\\
&      & +  \frac{1}{q}\left(1- \lambda_2 C_{m_2,q}^q(\Omega)\right)\|v\|_{W_0^{m_2,q}(\Omega)}^q
         -  \frac{\lambda_2(q-\gamma_2)}{q\gamma_2}\|h_2\|_{L^{\frac{q}{q-\gamma_2}}(\Omega)}^{\frac{q}{q-\gamma_2}}\nonumber\\
&      & -  \frac{C_0}{\alpha+\beta} \left(\frac{\alpha C_{m_1,p}^{\alpha+\beta}(\Omega)}{\alpha+\beta} \|u\|_{W_0^{m_1,p}(\Omega)}^{\alpha+\beta}
         +  \frac{\beta C_{m_2,q}^{\alpha+\beta}(\Omega)}{\alpha+\beta}  \|v\|_{W_0^{m_2,q}(\Omega)}^{\alpha+\beta}\right)\nonumber\\
& \ge  & \min\left\{ \frac{1-\lambda_1 C_{m_1,p}^p(\Omega)}{p} ,\frac{1- \lambda_2C_{m_2,q}^q(\Omega)}{q}\right\}
          \left(\|u\|_{W_0^{m_1,p}(\Omega)}^{\max\{p,q\}}  +  \|v\|_{W_0^{m_2,q}(\Omega)}^{\max\{p,q\}}\right)\nonumber\\
&      & - \frac{C_0}{(\alpha+\beta)^2} \left(\alpha C_{m_1,p}^{\alpha+\beta}(\Omega) + \beta C_{m_2,q}^{\alpha+\beta}(\Omega) \right)\|(u,v)\|_W^{\alpha+ \beta}\nonumber\\
&      &  -  \frac{\lambda_1(p-\gamma_1)} {p\gamma_1}\|h_1\|_{L^{\frac{p}{p-\gamma_1}}(\Omega)}^{\frac{p}{p-\gamma_1}}
           -  \frac{\lambda_2(q-\gamma_2)}{q\gamma_2}\|h_2\|_{L^{\frac{q}{q-\gamma_2}}(\Omega)}^{\frac{q}{q-\gamma_2}}   \nonumber\\
& \ge  & 2^{1-\max\{p,q\}}\min\left\{ \frac{1-\lambda_1 C_{m_1,p}^p(\Omega)}{p} ,\frac{1- \lambda_2C_{m_2,q}^q(\Omega)}{q}\right\}\|(u,v)\|_W^{\max\{p,q\}}\nonumber\\
&      & - \frac{C_0}{(\alpha+\beta)^2} \left(\alpha C_{m_1,p}^{\alpha+\beta}(\Omega) + \beta C_{m_2,q}^{\alpha+\beta}(\Omega) \right)\|(u,v)\|_W^{\alpha+ \beta}\nonumber\\
&      & -  \frac{\lambda_1(p-\gamma_1)} {p\gamma_1}\|h_1\|_{L^{\frac{p}{p-\gamma_1}}(\Omega)}^{\frac{p}{p-\gamma_1}}
         -   \frac{\lambda_2(q-\gamma_2)}{q\gamma_2}\|h_2\|_{L^{\frac{q}{q-\gamma_2}}(\Omega)}^{\frac{q}{q-\gamma_2}} .
\end{eqnarray}
Note that
\begin{eqnarray*}
M_{(\lambda_1,\lambda_2)} & = & 2^{1-\max\{p,q\}}\min\left\{ \frac{1-\lambda_1 C_{m_1,p}^p(\Omega)}{p},\frac{1- \lambda_2C_{m_2,q}^q(\Omega)}{q}\right\},\\
M_2 & = & \frac{C_0}{(\alpha+\beta)^2} \left(\alpha C_{m_1,p}^{\alpha+\beta}(\Omega) + \beta C_{m_2,q}^{\alpha+\beta}(\Omega) \right).
\end{eqnarray*}
Define
\begin{eqnarray}
\label{eq35}
f(t)=    M_{(\lambda_1,\lambda_2)} t^{\max\{p,q\}}- M_2 t^{\alpha+\beta}
     -  \frac{\lambda_1(p-\gamma_1)} {p\gamma_1}\|h_1\|_{L^{\frac{p}{p-\gamma_1}}(\Omega)}^{\frac{p}{p-\gamma_1}}
     -  \frac{\lambda_2(q-\gamma_2)}{q\gamma_2}\|h_2\|_{L^{\frac{q}{q-\gamma_2}}(\Omega)}^{\frac{q}{q-\gamma_2}},\ \ t \in [0,1].
\end{eqnarray}
To find $\rho_{(\lambda_1,\lambda_2)}$ such that $\varphi(u,v)>0$ whenever $\|(u,v)\|_W=\rho_{(\lambda_1,\lambda_2)}$, it is suffices to show that
there exists $t_{(\lambda_1,\lambda_2)}^{\star}\in (0,1]$ such that $f(t_{(\lambda_1,\lambda_2)}^{\star})>0$.
In fact, by (\ref{eq35}), we have
\begin{eqnarray*}
f'(t) & = & \max\{p,q\} M_{(\lambda_1,\lambda_2)} t^{\max\{p,q\}-1} - (\alpha+\beta) M_2 t^{\alpha+\beta-1},\\
f''(t) & = & \max\{p,q\}(\max\{p,q\}-1) M_{(\lambda_1,\lambda_2)} t^{\max\{p,q\}-2} - (\alpha+\beta)(\alpha+\beta-1) M_2 t^{\alpha+\beta-2}.
\end{eqnarray*}
Let $f'(t_{(\lambda_1,\lambda_2)}^{\star})=0$. We obtain that
$$
t_{(\lambda_1,\lambda_2)}^{\star}= \left(\frac{\max\{p,q\} M_{(\lambda_1,\lambda_2)}}{(\alpha+\beta)M_2}\right)^{\frac{1}{\alpha+\beta-\max\{p,q\}}}.
$$
Since $\lambda_1$ and $\lambda_2$ satisfy (\ref{eqq0}), we have
$
0<t_{(\lambda_1,\lambda_2)}^{\star} \le 1.
$
Moreover,
\begin{eqnarray*}
f''(t_{(\lambda_1,\lambda_2)}^{\star})
& = &   \max\{p,q\}(\max\{p,q\}-1) M_{(\lambda_1,\lambda_2)} \left(\frac{\max\{p,q\} M_{(\lambda_1,\lambda_2)}}{(\alpha+\beta)M_2}\right)^{\max\{p,q\}-2} \\
&   &  - (\alpha+\beta)(\alpha+\beta-1) M_2  \left(\frac{\max\{p,q\} M_{(\lambda_1,\lambda_2)}}{(\alpha+\beta)M_2}\right)^{\alpha+\beta-2}\\
& = &   (\max\{p,q\}-\alpha-\beta) M_2  \left(\frac{\max\{p,q\} M_{(\lambda_1,\lambda_2)}}{(\alpha+\beta)M_2}\right)^{\alpha+\beta-2}\\
& < & 0.
\end{eqnarray*}
Hence, by (\ref{eqq0}), we have
\begin{eqnarray*}
       \max_{t\in [0,1]}f(t)
& = & f(t_{(\lambda_1,\lambda_2)}^{\star})\\
&  =  &     M_{(\lambda_1,\lambda_2)} \left(\frac{\max\{p,q\} M_{(\lambda_1,\lambda_2)}}{(\alpha+\beta)M_2}\right)^{\frac{\max\{p,q\}}{\alpha+\beta-\max\{p,q\}}}
         -  M_2 \left(\frac{\max\{p,q\} M_{(\lambda_1,\lambda_2)}}{(\alpha+\beta)M_2}\right)^{\frac{\alpha+\beta}{\alpha+\beta-\max\{p,q\}}} \\
&     &  -  \frac{\lambda_1(p-\gamma_1)} {p\gamma_1}\|h_1\|_{L^{\frac{p}{p-\gamma_1}}(\Omega)}^{\frac{p}{p-\gamma_1}}
         -  \frac{\lambda_2(q-\gamma_2)}{q\gamma_2}\|h_2\|_{L^{\frac{q}{q-\gamma_2}}(\Omega)}^{\frac{q}{q-\gamma_2}}\\\
&  =  & \frac{\alpha+\beta - \max\{p,q\}}{\alpha+\beta}
        M_{(\lambda_1,\lambda_2)}^{\frac{\alpha+\beta}{\alpha+\beta-\max\{p,q\}}}
        \left(\frac{\max\{p,q\}}{(\alpha+\beta)M_2}\right)^{\frac{\max\{p,q\}}{\alpha+\beta-\max\{p,q\}}} \\
&     & -  \frac{\lambda_1(p-\gamma_1)} {p\gamma_1}\|h_1\|_{L^{\frac{p}{p-\gamma_1}}(\Omega)}^{\frac{p}{p-\gamma_1}}
        -   \frac{\lambda_2(q-\gamma_2)}{q\gamma_2}\|h_2\|_{L^{\frac{q}{q-\gamma_2}}(\Omega)}^{\frac{q}{q-\gamma_2}}\\
&  >  & 0.
\end{eqnarray*}
Let $\rho_{(\lambda_1,\lambda_2)}=t_{(\lambda_1,\lambda_2)}^{\star}$.
Thus, we conclude that  $\varphi(u,v)>0$ when $\|(u,v)\|_W= \rho_{(\lambda_1,\lambda_2)}$. \qed

\vskip2mm
\noindent
{\bf Lemma 3.2.} {\it For each $(\lambda_1,\lambda_2)$ satisfying (\ref{eqq0}),
there exists a $(u_{(\lambda_1,\lambda_2)}^{\star},v_{(\lambda_1,\lambda_2)}^{\star})\in W$ with $\|(u_{(\lambda_1,\lambda_2)}^{\star},v_{(\lambda_1,\lambda_2)}^{\star})\|_W>\rho_{(\lambda_1,\lambda_2)}$ such that $\varphi(u_{(\lambda_1,\lambda_2)}^{\star},v_{(\lambda_1,\lambda_2)}^{\star})<0$.}\\
{\bf Proof.} For any given $(u,v)\in W$ with $\int_{\Omega} c(x)|u|^{\alpha}|v|^{\beta}d\mu\not=0$ and any $z\in \R^+$, we have
\begin{eqnarray}
\label{eq321}
           \varphi(zu,zv)
&  =   &    \frac{1}{p}z^p \|u\|_{W_0^{m_1,p}(\Omega)}^p - \frac{\lambda_1}{\gamma_1}z^{\gamma_1}\int_{\Omega}h_1(x)|u|^{\gamma_1}d\mu\nonumber\\
&      &  + \frac{1}{q}z^q \|v\|_{W_0^{m_2,q}(\Omega)}^q - \frac{\lambda_2}{\gamma_2}z^{\gamma_2}\int_{\Omega}h_2(x)|v|^{\gamma_2}d\mu - \frac{1}{\alpha+\beta}z^{\alpha+\beta}\int_{\Omega} c(x)|u|^{\alpha}|v|^{\beta}d\mu\nonumber\\
& \le  &    \frac{1}{p}z^p \|u\|_{W_0^{m_1,p}(\Omega)}^p - \frac{\lambda_1h_1^{\star}}{\gamma_1}z^{\gamma_1}\int_{\Omega}|u|^{\gamma_1}d\mu\nonumber\\
&      &  + \frac{1}{q}z^q \|v\|_{W_0^{m_2,q}(\Omega)}^q - \frac{\lambda_2h_2^{\star}}{\gamma_2}z^{\gamma_2}\int_{\Omega}|v|^{\gamma_2}d\mu - \frac{1}{\alpha+\beta}z^{\alpha+\beta}\int_{\Omega} c(x)|u|^{\alpha}|v|^{\beta}d\mu.
\end{eqnarray}
Note that $\alpha+\beta>\max\{p,q\}$. So, there exists $z_{(\lambda_1,\lambda_2)}^{\star}$ large enough such that $\|(z_{(\lambda_1,\lambda_2)}^{\star}u,z_{(\lambda_1,\lambda_2)}^{\star}v)\|_W>\rho_{(\lambda_1,\lambda_2)}$ and $\varphi(z_{(\lambda_1,\lambda_2)}^{\star} u,z_{(\lambda_1,\lambda_2)}^{\star} v)<0$.
Let $u_{(\lambda_1,\lambda_2)}^{\star}=z_{(\lambda_1,\lambda_2)}^{\star} u$ and $v_{(\lambda_1,\lambda_2)}^{\star}=z_{(\lambda_1,\lambda_2)}^{\star} v$. Then the proof is finished.
\qed

\vskip2mm
\noindent
{\bf Lemma 3.3.} {\it For each $(\lambda_1,\lambda_2)$ satisfying (\ref{eqq0}), $\varphi$ satisfies the  Palais-Smale condition, namely, any sequence $\{(u_k,v_k)\}\subset W $, which satisfies that $\varphi(u_k,v_k)$ is bounded for all $k \in \mathbb{N}$ and $ \varphi'(u_k,v_k)\rightarrow 0$ as $ k\to\infty$, has a convergent subsequence. } \\
{\bf Proof.} Assume that  there exists a positive constant $c$ such that
$$
|\varphi(u_k,v_k)|\le c \ \ \mbox{and} \ \ \varphi'(u_k,v_k)\rightarrow 0 \ \ \mbox{as } k\to\infty, \ \ \mbox{ for all }k\in \mathbb{N}  .
$$
Then
\begin{eqnarray}
\label{eq331}
&       &   c+\|u_k\|_{W_0^{m_1,p}(\Omega)}+\|v_k\|_{W_0^{m_2,q}(\Omega)}\nonumber\\
&  \ge  &   \varphi(u_k,v_k)-\frac{1}{\alpha+\beta}\langle\varphi'(u_k,v_k),(u_k,v_k)\rangle\nonumber\\
&   =   &   \left(\frac{1}{p} - \frac{1}{\alpha+\beta}\right)\|u_k\|_{W_0^{m_1,p}(\Omega)}^p + \left(\frac{1}{q}-\frac{1}{\alpha+\beta}\right)\|v_k\|_{W_0^{m_2,q}(\Omega)}^q \nonumber\\
&       &  - \lambda_1\left(\frac{1}{\gamma_1}-\frac{1}{\alpha+\beta}\right)\int_{\Omega}h_1(x)|u_k|^{\gamma_1}d\mu - \lambda_2\left(\frac{1}{\gamma_2}-\frac{1}{\alpha+\beta}\right)\int_{\Omega}h_2(x)|v_k|^{\gamma_2}d\mu\nonumber\\
&  \ge  &   \left(\frac{1}{p} - \frac{1}{\alpha+\beta}\right)\|u_k\|_{W_0^{m_1,p}(\Omega)}^p + \left(\frac{1}{q}-\frac{1}{\alpha+\beta}\right)\|v_k\|_{W_0^{m_2,q}(\Omega)}^q \nonumber\\
&       &  - \lambda_1\left(\frac{1}{\gamma_1}-\frac{1}{\alpha+\beta}\right)H_1C_{m_1,p}^{\gamma_1}(\Omega)\|u_k\|_{W_0^{m_1,p}(\Omega)}^{\gamma_1} \nonumber\\
& &- \lambda_2\left(\frac{1}{\gamma_2}-\frac{1}{\alpha+\beta}\right)H_2C_{m_2,q}^{\gamma_2}(\Omega)\|v_k\|_{W_0^{m_2,q}(\Omega)}^{\gamma_2}.
\end{eqnarray}
We claim that ${\|(u_k,v_k)\|}_W$ is bounded. In fact,  if
\begin{eqnarray}
\label{eq332}
\|u_k\|_{W_0^{m_1,p}(\Omega)}\to\infty  \mbox{ and } \|v_k\|_{W_0^{m_2,q}(\Omega)}\to \infty \mbox{ as } k\to \infty,
\end{eqnarray}
then it follows from (\ref{eq331})  that
\begin{eqnarray}
\label{eq333}
&       &   c+\|(u_k,v_k)\|_W\nonumber\\
&  \ge  &   \min\left\{\frac{1}{p} - \frac{1}{\alpha+\beta}, \frac{1}{q}-\frac{1}{\alpha+\beta}\right\}\left(\|u_k\|_{W_0^{m_1,p}(\Omega)}^p +\|v_k\|_{W_0^{m_2,q}(\Omega)}^q \right)\nonumber\\
&       &  - \max\left\{\lambda_1\left(\frac{1}{\gamma_1}-\frac{1}{\alpha+\beta}\right)H_1C_{m_1,p}^{\gamma_1}(\Omega),\lambda_2\left( \frac{1}{\gamma_2} - \frac{1}{\alpha+\beta}\right)H_2C_{m_2,q}^{\gamma_2}(\Omega)\right\}\left(\|u_k\|_{W_0^{m_1,p}(\Omega)}^{\gamma_1}+\|v_k\|_{W_0^{m_2,q}(\Omega)}^{\gamma_2}\right)\nonumber\\
&  \ge  &    2^{1-\min\{p,q\}}\min\left\{\frac{1}{p} - \frac{1}{\alpha+\beta}, \frac{1}{q}-\frac{1}{\alpha+\beta}\right\}\|(u_k,v_k)\|_W^{\min\{p,q\}} \nonumber\\
&       &  - \max\left\{\lambda_1\left(\frac{1}{\gamma_1}-\frac{1}{\alpha+\beta}\right)H_1C_{m_1,p}^{\gamma_1}(\Omega),\lambda_2\left( \frac{1}{\gamma_2} - \frac{1}{\alpha+\beta}\right)H_2C_{m_2,q}^{\gamma_2}(\Omega)\right\}\|(u_k,v_k)\|_W^{\max\{\gamma_1,\gamma_2\}},\nonumber
\end{eqnarray}
which contradicts with (\ref{eq332}). If
\begin{eqnarray}
\label{eq334}
\|u_k\|_{W_0^{m_1,p}(\Omega)}\to\infty \mbox{ as } k\to \infty,
\end{eqnarray}
and $\|v_k\|_{W_0^{m_2,q}(\Omega)}$ is bounded for all $k\in \mathbb{N}$, then it follows from (\ref{eq331}) that there exists a positive constant $c_1$ such that
\begin{eqnarray*}
\label{eq335}
      c_1+\|u_k\|_{W_0^{m_1,p}(\Omega)}
 \ge  \left(\frac{1}{p} - \frac{1}{\alpha+\beta}\right)\|u_k\|_{W_0^{m_1,p}(\Omega)}^p
      -\lambda_1\left(\frac{1}{\gamma_1}-\frac{1}{\alpha+\beta}\right)H_1C_{m_1,p}^{\gamma_1}(\Omega)\|u_k\|_{W_0^{m_1,p}(\Omega)}^{\gamma_1} ,
\end{eqnarray*}
which contradicts with (\ref{eq334}). Similarly, if
\begin{eqnarray*}
\label{eq336}
\|v_k\|_{W_0^{m_2,q}(\Omega)}\to\infty \mbox{ as } k\to \infty,
\end{eqnarray*}
and $\|u_k\|_{W_0^{m_1,p}(\Omega)}$ is bounded for all $k\in \mathbb{N}$, we can also obtain a contradiction.
Hence, $\|u_k\|_{W_0^{m_1,p}(\Omega)}$ and  $\|v_k\|_{W_0^{m_2,q}(\Omega)}$ are bounded for all $k \in \mathbb{N}$.
Then, there is a subsequence, still denoted by $\{u_k\}$, such that $u_k\rightharpoonup u^{\star}$ for some $u^{\star}\in W_0^{m_1,p}(\Omega)$ as $k\to \infty$,
and a subsequence of $\{v_k\}$,  still denoted by $\{v_k\}$,
such that $v_k\rightharpoonup v^{\star}$ for some $v^{\star}\in W_0^{m_2,q}(\Omega)$ as $k\to \infty$.
 Lemma 2.1  implies that
\begin{eqnarray*}
\label{eq338}
 u_k\to  u^{\star} \ \mbox{in } W_0^{m_1,p}(\Omega)\ \ \mbox{and} \ \ v_k\to  v^{\star} \ \mbox{in } W_0^{m_2,q}(\Omega)\ \ \mbox{as }k\to \infty.
\end{eqnarray*}
The proof is complete.
\qed

\vskip2mm
 \noindent
 {\bf Proof of Theorem 1.1.} By Lemma 3.1-Lemma 3.3 and Lemma 2.3, we obtain that for each $(\lambda_1,\lambda_2)$ satisfying (\ref{eqq0}), system (\ref{eq1}) has one nontrivial solution $(u^{\star},v^{\star})$  of positive energy.
 Next, we prove that system (\ref{eq1}) has one nontrivial solution of  negative energy. The proof is similar to Theorem 3.3 in \cite{Cheng} and Theorem 1.3 in \cite{Yang-Zhang}. Note that $\gamma_1<p$ and $\gamma_2<q$. Then it follows from (\ref{eq321}) that there exists $z$ small enough such that
$$
 \varphi(zu,zv)<0.
$$
So,
$$
 -\infty<\inf\{\varphi(u,v):(u,v)\in\bar B_{\rho_{(\lambda_1,\lambda_2)}} \}<0,
$$
where $\rho_{(\lambda_1,\lambda_2)}$ is given in Lemma 3.1 and $\bar B_{\rho_{(\lambda_1,\lambda_2)}}=\{(u,v)\in W\big|\|(u,v)\|_W\le \rho_{(\lambda_1,\lambda_2)}\}$.
Moreover, by Lemma 3.1, we have
$$
-\infty<\inf_{\bar B_{\rho_{(\lambda_1,\lambda_2)}}}\varphi(u,v)<0<\inf_{\partial B_{\rho_{(\lambda_1,\lambda_2)}}}\varphi(u,v)
$$
for each $(\lambda_1,\lambda_2)$ satisfying (\ref{eqq0}). Set
\begin{eqnarray}
\label{A1}
\frac{1}{n}\in \left(0,\inf_{\partial B_{\rho_{(\lambda_1,\lambda_2)}}}\varphi(u,v)-\inf_{\bar B_{\rho_{(\lambda_1,\lambda_2)}}} \varphi(u,v)\right),\; n\in \mathbb{N}.
\end{eqnarray}
It follows from the definition of infimum that there exists a $(u_n,v_n)\in \bar B_{\rho_{(\lambda_1,\lambda_2)}}$ such that
\begin{eqnarray}
\label{EQ11}
\varphi(u_n,v_n)\le\inf_{\bar B_{\rho_{(\lambda_1,\lambda_2)}}}\varphi(u,v)+\frac {1}{n}.
\end{eqnarray}
Note that $\varphi(u,v) \in C^1(W,\R)$ and then $\varphi(u,v)$ is lower semicontinuous. Hence, using Lemma 2.4, we get
$$
\varphi(u_n,v_n)\le \varphi(u,v)+\frac{1}{n}\|(u,v)-(u_n,v_n)\|_W,\;\forall (u,v)\in \bar B_{\rho_{(\lambda_1,\lambda_2)}}.
$$
(\ref{A1}) and (\ref{EQ11}) imply that
$$
\varphi(u_n,v_n) \le \inf_{\bar B_{\rho_{(\lambda_1,\lambda_2)}}}\varphi(u,v)+\frac {1}{n}<\inf_{\partial B_{\rho_{(\lambda_1,\lambda_2)}}}\varphi(u,v),
$$
 and so $(u_n,v_n)\in B_{\rho_{(\lambda_1,\lambda_2)}}$. Set $ M_n:W \rightarrow \R$ as
$$
M_n(u,v)=\varphi(u,v) + \frac {1}{n}\|(u,v)-(u_n,v_n)\|_W.
$$
Then $(u_n,v_n)\in B_{\rho_{(\lambda_1,\lambda_2)}}$ is the minimum point of $M_n$ on $ \bar B_{\rho_{(\lambda_1,\lambda_2)}}$. Hence, for some $(u,v)\in W$ satisfying $\|(u,v)\|_W = 1$,
if $t >0$ is small enough such that $(u_n+tu,v_n+tv)\in \bar B_{\rho_{(\lambda_1,\lambda_2)}}$, then
\begin{eqnarray}
\label{A2}
\frac{M_n(u_n+tu,v_n+tv)-M_n(u_n,v_n)}{t}\ge 0.
\end{eqnarray}
By (\ref{A2}) and the definition of $M_n$, there holds
$$
\langle\varphi'(u_n,v_n),(u,v)\rangle\ge-\frac{1}{n}.
$$
Similarly, if $t<0$ and $|t|$ is small enough, then
$$
\langle\varphi'(u_n,v_n),(u,v)\rangle\le\frac{1}{n}.
$$
Therefore,
\begin{eqnarray}
\label{EQ12}
\|\varphi'(u_n,v_n)\|=\sup_{\|(u,v)\|_W=1}|\langle\varphi'(u_n,v_n),(u,v)\rangle|\le\frac{1}{n}.
\end{eqnarray}
Thus (\ref{EQ11}) and (\ref{EQ12}) imply  that
$$
\varphi(u_n,v_n)\rightarrow \inf_{\bar B_{\rho_{(\lambda_1,\lambda_2)}}} \varphi(u,v)\;\;\;\mbox{and} \;\;\;\|\varphi'(u_n,v_n)\|\rightarrow 0\;\;\;as \;n\rightarrow\infty.
$$
Then it follows from Lemma 3.3 that $\{(u_n,v_n)\}$ has a strongly convergent subsequence, still denoted by $(u_{n},v_{n})$,
 and $(u_{n},v_{n})\rightarrow(u^{\star\star},v^{\star\star})\in \bar B_{\rho_{(\lambda_1,\lambda_2)}}$ as $n \rightarrow\infty$, and
$$
\varphi(u^{\star\star},v^{\star\star})=\inf_{\bar B_{\rho_{(\lambda_1,\lambda_2)}}}\varphi<0\;\;\;\mbox{and} \;\;\;\varphi'(u^{\star\star},v^{\star\star})=0.
$$
Hence, system (\ref{eq1}) has a nontrivial solution $(u^{\star\star},v^{\star\star})$ of  negative energy.
The proof is complete. \qed

\vskip2mm
\noindent
{\bf Proof of Theorem 1.2.} For each $\lambda_1>0$,  if  $(u,v)=(u,0)$ is a semi-trivial solution of system (\ref{eq1}), then we have
\begin{eqnarray*}
\label{EQ13}
      \int_{\Omega \cup \partial\Omega}|\nabla^{m_1} u|^p d\mu
  =   \lambda_1\int_{\Omega}h_1(x)|u|^{\gamma_1}d\mu
 \le  \lambda_1 H_1 \int_{\Omega}|u|^{\gamma_1}d\mu
 \le  \lambda_1 H_1 C_{m_1,p}^{\gamma_1}(\Omega)\|u\|_{W_0^{m_1,p}(\Omega)}^{\gamma_1}.
\end{eqnarray*}
Hence,
 $$
 \|u\|_{W_0^{m_1,p}(\Omega)} \le \left(\lambda_1 H_1C^{\gamma_1}_{m_1,p}(\Omega)\right)^{\frac{1}{p-\gamma_1}}.
 $$
Similarly, for each $\lambda_2>0$,  if $(u,v)=(0,v)$  is a semi-trivial solution of system (\ref{eq1}), we can also obtain
 $$
 \|v\|_{W_0^{m_2,q}(\Omega)} \le \left(\lambda_2 H_2C^{\gamma_2}_{m_2,q}(\Omega)\right)^{\frac{1}{q-\gamma_2}}.
 $$
The proof is complete.
\qed

\section{Proofs for  Theorem 1.3}
\setcounter{equation}{0}
In this section, we discuss the existence of  ground state solution for equation (\ref{a1}) by using Lemma 2.5 and Theorem 3.3 in \cite{Brown 09}. In \cite{Brown 09}, Brown and Wu researched the following operator equation with the help of fibering maps and the Nehari manifold:
\begin{eqnarray}
\label{500}
 A(u)-B(u)-C(u)=0, \ \ u\in X,
\end{eqnarray}
where $X$ is a reflexive Banach space, $A,B,C: X \to X^*$ are homogeneous operators of degree $p-1,\alpha-1$ and $\gamma-1$ with $1<\gamma<p<\alpha$.
The energy functional of (\ref{500}) is
\begin{eqnarray}
\label{eq521}
J(u) =   \frac{1}{p}         \langle A(u), u\rangle
         - \frac{1}{\alpha}  \langle B(u), u\rangle
         - \frac{1}{\gamma}  \langle C(u), u\rangle,
\end{eqnarray}
the fibering map is
\begin{eqnarray}
\label{eq522}
G_u(t) =    \frac{1}{p}       \langle A(u), u\rangle t^p
         - \frac{1}{\alpha}  \langle B(u), u\rangle t^{\alpha}
         - \frac{1}{\gamma } \langle C(u), u\rangle t^{\gamma},
\end{eqnarray}
for all $t>0$, and the Nehari manifold is defined by
$$
\mathcal{N}=\{u \in X \setminus \{0\}\big|  \langle J'(u),u\rangle=0\},
$$
Define
$$
\phi(u)=\langle J'(u),u\rangle.
$$
Then $\mathcal{N}$  can be divided into three parts:
\begin{eqnarray*}
\mathcal{N}^+=\{u\in \mathcal{N} \big| \langle \phi'(u),u\rangle>0\},\\
\mathcal{N}^0=\{u\in \mathcal{N} \big| \langle \phi'(u),u\rangle=0\},\\
\mathcal{N}^-=\{u\in \mathcal{N} \big| \langle \phi'(u),u\rangle<0\}.
\end{eqnarray*}
The following results were obtained in \cite{Brown 09}.
\vskip2mm
\noindent
{\bf Theorem A.} (\cite{Brown 09}, Lemma 2.5) {\it For each $u \in X$, when $\langle B(u), u\rangle>0$ and $\langle C(u), u\rangle >0$,
there are $t_u^+$ and $t_u^-$ with $0< t_u^+< t_u^-$ such that
$G_u(t)$ is decreasing on $(0,t_u^+)$, increasing on $(t_u^+,t_u^-)$ and decreasing on $(t_u^-, +\infty)$.
}
\vskip0mm
\noindent
{\bf Theorem B.} (\cite{Brown 09}, Theorem 3.3) {\it Assume that the following conditions hold:\\
$(H_1)$ $u \to \langle A(u),u\rangle $  is weakly lower semicontinuous on $X$ and there exists a continuous function $\kappa: [0,+\infty) \to [0,+\infty)$ with $\kappa(s) >0$
 on $(0, +\infty)$ and $\lim_{s\rightarrow \infty}\kappa(s) = \infty$ such that $\langle A(u),u\rangle \ge \kappa(\|u\|)\|u\|$ for all $u \in X$;\\
$(H_2)$ there exist $u_i \in X$, $i=1,2$  such that
 $$
 \langle B(u_1), u_1\rangle>0,
 \langle C(u_2), u_2\rangle>0,
 $$
$(H_3)$ $B,C$ are strongly continuous;\\
$(H_4)$ there exist two positive constants $d_1,d_2$ with
$$
d_1^{\alpha-p} d_2^{p-\gamma}  \le (p-\gamma)^{p-\gamma}(\alpha-p)^{\alpha-p}(\alpha-\gamma)^{\gamma-\alpha},
$$
such that
\begin{eqnarray}
\label{501}
   \langle B(u), u\rangle \le d_1[\langle A(u),u \rangle]^{\frac{\gamma}{p}},\\
\label{502}
 \langle C(u_), u\rangle \le d_2[\langle A(u),u \rangle]^{\frac{\alpha}{p}},
 \end{eqnarray}
equation (\ref{500}) has at least two nontrivial solutions $u_0^+$ and $u_0^-$, where
\begin{eqnarray*}
u_0^+ \in \mathcal{N}^+,\ \ J(u_0^+)=\inf_{u\in \mathcal{N}^+}J(u),\\
u_0^- \in \mathcal{N}^-,\ \ J(u_0^-)=\inf_{u\in \mathcal{N}^-}J(u),
\end{eqnarray*}
and $\mathcal{N}^0=\emptyset$.
}
\vskip2mm
\par
In the locally finite graph $G=(V,E)$ setting, let $X=W_0^{m,p}(\Omega)$ and
\begin{eqnarray}
\label{1}\langle A(u), u \rangle & = & \|u\|_{W_0^{m,p}(\Omega)}^p,\\
\label{2}\langle B(u), u \rangle & = & \int_{\Omega}c(x)|u|^{\alpha}d\mu,\\
\label{3}\langle C(u), u \rangle & = & \lambda\int_{\Omega}h(x)|u|^{\gamma}d\mu.
\end{eqnarray}
Similar to the arguments in \cite{Brown 09}, $(H_1)$-$(H_3)$ hold with $A,B,C$ defined by (\ref{1})-(\ref{3}).
Note that
$$
C_0=\max_{x \in \Omega} c(x),\ \  H_0=\max_{x \in \Omega} h(x).
$$
By Lemma 2.1, we have
\begin{eqnarray}
\label{eq41}
     \int_{\Omega}c(x)|u|^{\alpha}d\mu
\le  C_0 C_{m,p}^{\alpha}(\Omega) \|u\|_{W_0^{m,p}(\Omega)}^{\alpha}
\end{eqnarray}
and
\begin{eqnarray}
\label{eq42}
       \lambda\int_{\Omega}h(x)|u|^{\gamma}d\mu
 \le   \lambda H_0 C_{m,p}^{\gamma}(\Omega) \|u\|_{W_0^{m,p}(\Omega)}^{\gamma}.
\end{eqnarray}
Let
\begin{eqnarray*}
d_1  =   C_0 C_{m,p}^{\alpha}(\Omega),\ \
d_2  =   \lambda H_0 C_{m,p}^{\gamma}(\Omega).
\end{eqnarray*}
Then (\ref{501}) and (\ref{502}) hold.
Moreover, note that
$$
\lambda_0
 =  \frac{p-\gamma}{H_0} C_{m,p}^{-\gamma}(\Omega) \left(\left(C_0 C_{m,p}^{\alpha}(\Omega) \right)^{p-\alpha}(\alpha-p)^{\alpha-p}(\alpha-\gamma)^{\gamma-\alpha}\right)^{\frac{1}{p-\gamma}}.
$$
Then it is  easy to see that $(H_4)$ holds if $\lambda \in (0,\lambda_0)$.
Thus, by Theorem B, equation (\ref{a1}) has at least two nontrivial solutions $u^+_0\in \mathcal{N}^+$ and $u^-_0\in\mathcal{N}^-$,
and one of  $u^+_0$ and $u^-_0$ must be ground state solution of equation (\ref{a1}).
Next, we discuss which  is the ground state solution.
Note that the energy functional of equation (\ref{a1}) is
\begin{eqnarray}
\label{eq522}
J(u) & = &   \frac{1}{p} \|u\|_{W_0^{m,p}(\Omega)}^p - \frac{\lambda}{\gamma}\int_{\Omega}h(x)|u|^{\gamma}d\mu - \frac{1}{\alpha}\int_{\Omega}c(x)|u|^{\alpha}d\mu, \ \ \forall u\in W_0^{m,p}(\Omega),
\end{eqnarray}
and for each  $u\in W_0^{m,p}(\Omega)\backslash\{0\}$, the corresponding fibering map is
\begin{eqnarray}
\label{eq523}
G_u(t) =    \frac{t^p }{p}     \|u\|_{W_0^{m,p}(\Omega)}^p
         - \frac{\lambda}{\gamma}t^{\gamma}  \int_{\Omega}h(x)|u|^{\gamma}d\mu
         - \frac{t^{\alpha}}{\alpha}   \int_{\Omega}c(x)|u|^{\alpha}d\mu,\ \   \forall t\in (0,+\infty).
\end{eqnarray}
We can obtain that $G_u(t)$ has positive values if $\lambda \in (0,\lambda_{\star})$ where $\lambda_{\star}$ is defined by (\ref{aa2}).
In fact, we define
$$
F_u(t) =   \frac{t^p}{p}    \|u\|_{W_0^{m,p}(\Omega)}^p
         - \frac{ t^{\alpha}}{\alpha}  \int_{\Omega}c(x)|u|^{\alpha}d\mu .
$$
By (\ref{eq41}), we have
\begin{eqnarray}
\label{eq531}
         \max_{t>0}F_{u}(t)
 &  =  & F_u(t_{0u})\\
 &  =  &  \frac{1}{p}    \|u\|_{W_0^{m,p}(\Omega)}^p \left(\frac{\|u\|_{W_0^{m,p}(\Omega)}^p}{\int_{\Omega}c(x)|u|^{\alpha}d\mu}\right)^{\frac{p}{\alpha-p}}
       - \frac{1}{\alpha}  \int_{\Omega}c(x)|u|^{\alpha}d\mu  \left(\frac{\|u\|_{W_0^{m,p}(\Omega)}^p}{\int_{\Omega}c(x)|u|^{\alpha}d\mu}\right)^{\frac{\alpha}{\alpha-p}}\nonumber\\
 &  =  &  \left(\frac{1}{p} -  \frac{1}{\alpha}\right) \left(\frac{\|u\|_{W_0^{m,p}(\Omega)}^{\alpha}}{\int_{\Omega}c(x)|u|^{\alpha}d\mu}\right)^{\frac{p}{\alpha-p}}\nonumber\\
 & \ge &  \frac{\alpha-p }{p\alpha} \left(C_0C^{\alpha}_{m,p}(\Omega)\right)^{\frac{ p}{p-\alpha}}
\end{eqnarray}
with
$$
t_{0u}=\left(\frac{\|u\|_{W_0^{m,p}(\Omega)}^p}{\int_{\Omega}c(x)|u|^{\alpha}d\mu}\right)^{\frac{1}{\alpha-p}}.
$$
Furthermore,
\begin{eqnarray}
\label{eq532}
&     &  \frac{\lambda}{\gamma}  t_{0u}^{\gamma} \int_{\Omega}h(x)|u|^{\gamma}d\mu\nonumber\\
& \le &  \frac{\lambda H_0}{\gamma}  C_{m,p}^\gamma(\Omega)\|u\|_{W_0^{m,p}(\Omega)}^{\gamma}t_{0u}^{\gamma}\nonumber\\
&  =  &  \frac{\lambda H_0}{\gamma} C_{m,p}^\gamma(\Omega) \|u\|_{W_0^{m,p}(\Omega)}^{\gamma} \left(\frac{\|u\|_{W_0^{m,p}(\Omega)}^p}{\int_{\Omega} c(x)|u|^{\alpha}d\mu}\right)^{\frac{\gamma}{\alpha-p}}\nonumber\\
&  =  &  \frac{\lambda H_0}{\gamma} C_{m,p}^\gamma(\Omega) \left(\frac{\|u\|_{W_0^{m,p}(\Omega)}^{\alpha}}{\int_{\Omega}c(x)|u|^{\alpha}d\mu}\right)^{\frac{\gamma}{\alpha-p}}\nonumber\\
&  =  &  \frac{\lambda H_0}{\gamma} C_{m,p}^\gamma(\Omega)\left(\frac{p\alpha}{\alpha-p} F_u(t_{0u})\right)^{\frac{\gamma}{p}}.
\end{eqnarray}
It follows that
\begin{eqnarray}
\label{eq533}
G_u(t_{0u})
&  =  & F_u(t_{0u}) - \frac{\lambda}{\gamma} t_{0u}^{\gamma} \int_{\Omega}h(x)|u|^{\gamma}d\mu \nonumber\\
& \ge & F_u^{\frac{\gamma}{p}}(t_{0u}) \left(F_u^{\frac{p-\gamma}{p}}(t_{0u}) - \lambda\frac{H_0}{\gamma}C_{m,p}^\gamma(\Omega)\left(\frac{p\alpha}{\alpha-p}\right)^{\frac{\gamma}{p}}\right).
\end{eqnarray}
Note that
$$
\lambda_{\star}= \frac{\gamma(\alpha-p)}{p\alpha H_0 C_{m,p}^\gamma(\Omega)} \left(C_0C^{\alpha}_{m,p}(\Omega)\right)^{\frac{p- \gamma}{p-\alpha}}.
$$
Then it is easy to see that for all $u \in W_0^{m,p}(\Omega)\backslash\{0\}$, if $\lambda \in (0, \lambda_{\star})$, there holds
\begin{eqnarray}\label{eqq1}
G_u(t_{0u})>0.
\end{eqnarray}
Moreover, for all $u \in \mathcal{N}^-$, by the definition of $\mathcal{N}^-$, we have
$$
G'_u(1)=  \|u\|_{W_0^{m,p}(\Omega)}^p -  \int_{\Omega}h(x)|u|^{\gamma}d\mu -   \int_{\Omega}c(x)|u|^{\alpha}d\mu = \langle J'(u), u \rangle = 0,
$$
and
\begin{eqnarray*}
G''_u(1)
& = & (p-1) \|u\|_{W_0^{m,p}(\Omega)}^p - (\gamma-1) \int_{\Omega}h(x)|u|^{\gamma}d\mu -  (\alpha - 1)\int_{\Omega}c(x)|u|^{\alpha}d\mu \\
& = & p \|u\|_{W_0^{m,p}(\Omega)}^p - \gamma \int_{\Omega}h(x)|u|^{\gamma}d\mu -  \alpha \int_{\Omega}c(x)|u|^{\alpha}d\mu -
\left(\|u\|_{W_0^{m,p}(\Omega)}^p -  \int_{\Omega}h(x)|u|^{\gamma}d\mu -   \int_{\Omega}c(x)|u|^{\alpha}d\mu \right)\\
& = & \langle \phi'(u), u \rangle - \langle J'(u), u \rangle\\
& < & 0.
\end{eqnarray*}
Then $t_u=1$ is a local maximum point of $G_u(t)$.
Note that for any $u \in W_0^{m,p}(\Omega)\backslash\{0\}$, there are
\begin{eqnarray*}
\langle B(u), u \rangle & = & \int_{\Omega}c(x)|u|^{\alpha}d\mu>0,\\
\langle C(u), u \rangle & = & \lambda\int_{\Omega}h(x)|u|^{\gamma}d\mu>0.
\end{eqnarray*}
Hence, by Theorem A, for each  $u \in \mathcal{N}^-$, there exist $t_u^+$ and $t_u^-$ with $0< t_u^+< t_u^-$ such that
$G_u(t)$ is decreasing on $(0,t_u^+)$, increasing on $(t_u^+,t_u^-)$ and decreasing on $(t_u^-, +\infty)$, together with $G_u(0)=0$ and (\ref{eqq1}), which implies that both the local maximum point  $t_u=1$ and $t_{0u}$ belong to the interval $(t_u^+, +\infty)$. Thus, we have
$$
J(u)=G_u(1) \ge G_u(t_{0u})>0,\ \ \mbox{for each }  u \in \mathcal{N}^-.
$$
Similarly, for all $u \in \mathcal{N}^+$, we know that $t_u=1$ is a local minimum point of $G_u(t)$, which locates at  $(0,t_u^-)$ by Theorem A. Hence, we have
$$
J(u)=G_u(1)<G_u(0)=0, \ \ \mbox{for each }  u \in \mathcal{N}^+.
$$
Hence, we conclude that  equation (\ref{a1}) has a nontrivial ground state solution $u^+_0\in \mathcal{N}^+$ if $\lambda \in (0, \lambda_{\star\star})$ where $\lambda_{\star\star}=\min\{\lambda_0,\lambda_{\star}\}$.
\qed

\section{Some results on the finite graph}
\setcounter{equation}{0}
If  $G=(V,E)$ is a finite graph, using the similar arguments as Theorem 1.1-Theorem 1.3,  we can obtain the similar results for the following poly-Laplacian system:
\begin{eqnarray}
\label{a2}
 \begin{cases}
  \pounds_{m_1,p}u+a(x)|u|^{p-2}u=\lambda_1 h_1(x)|u|^{\gamma_1-2}u+\frac{\alpha}{\alpha+\beta}c(x)|u|^{\alpha-2}u|v|^{\beta},&\;\;\;\;\hfill x\in V,\\
   \pounds_{m_2,q}v+b(x)|v|^{q-2}v=\lambda_2 h_2(x)|v|^{\gamma_2-2}v+\frac{\beta}{\alpha+\beta}c(x)|u|^{\alpha}|v|^{\beta-2}v,&\;\;\;\;\hfill x\in V,\\
   \end{cases}
\end{eqnarray}
where $m_i,  i=1,2$ are positive integers, $p,q,\gamma_1,\gamma_2>1$, $\lambda_1,\lambda_2,\alpha,\beta>0$, $\max\{\gamma_1,\gamma_2\} < \min\{p,q\} \le \max\{p,q\} < \alpha+\beta$, $a,b,h_1,h_2, c:V \to \R^+$,
and the following  equation:
\begin{eqnarray}
\label{a3}
  \pounds_{m,p}u+a(x)|u|^{p-2}u=\lambda h(x)|u|^{\gamma-2}u+c(x)|u|^{\alpha-2}u,&\;\;\;\;\hfill x\in V,
\end{eqnarray}
where $m$ is a positive integer, $p,\gamma>1$, $\lambda,\alpha>0$, $\gamma < p < \alpha$, $a,h, c:V \to \R^+$.
For any given $s>1$, we define $W^{m,s}(V)$ as a set of all functions $u:V\to \R$
with the norm
\begin{eqnarray*}
\|\psi\|_{W^{m,s}(V)}=\left(\int_V(|\nabla^{m} \psi(x)|^s+h(x)|\psi(x)|^s)d\mu\right)^\frac{1}{s},
\end{eqnarray*}
and for any given $1\le r<+\infty$, we define $L^r(V)$ as a set of all functions $u:V\to \R$ with the norm
$$
\|u\|_{L^r(V)}=
\left(\int_V|u(x)|^r d \mu\right)^{\frac{1}{r}},
$$
For system (\ref{a2}), we work in the space $W(V)=W^{m_1,p}(V)\times W^{m_2,q}(V)$, and for equation (\ref{a3}), we work in the space $W^{m,p}(V)$. Both $W(V)$ and $W^{m,p}(V)$ are of finite dimension. See \cite{Yamabe 2016} for more details.

\par
Denote
\begin{eqnarray*}
M_{(\lambda_1,\lambda_2)}(V) & = & 2^{1-\max\{p,q\}}\min\left\{ \frac{1-\lambda_1 C_p^p(V)}{p},\frac{1- \lambda_2C_q^q(V)}{q}\right\},\\
M_2(V) & = & \frac{C_0(V)}{(\alpha+\beta)^2} \left(\alpha C_p^{\alpha+\beta}(V) + \beta C_q^{\alpha+\beta}(V) \right),
\end{eqnarray*}
where $C_0(V)=\max_{x \in V}c(x)$,  $C_p(V)$ and $C_q(V)$ are embedding constants from $W^{m,p}(V)$ and $W^{m,q}(V)$ into $L^p(V)$ and $L^q(V)$, respectively, which have been obtained in \cite{Zhang-Zhang} with
\begin{eqnarray*}
 C_p(V)=\frac{\left(\sum_{x\in V}\mu(x)\right)^{\frac{1}{p}}}{\mu_{\min}^{\frac{1}{p}}h_{\min}^{\frac{1}{p}}}
\mbox{ and }  C_q(V)=\frac{\left(\sum_{x\in V}\mu(x)\right)^{\frac{1}{q}}}{\mu_{\min}^{\frac{1}{q}}h_{\min}^{\frac{1}{q}}}.
\end{eqnarray*}
 \par
 Next, we state the results similar to Theorem 1.1-Theorem 1.3. Assume that $\lambda_1$ and $\lambda_2$ satisfy the following inequalities:
\begin{eqnarray}
\begin{cases}
\label{eqqq0}
 &  0 < \lambda_1 < C_p^{-p}(V),\\
 &  0 < \lambda_2 < C_q^{-q}(V),\\
 &  M_{(\lambda_1,\lambda_2)}(V) \le \frac{\alpha+\beta}{\max\{p,q\}}M_2,\\
 & \frac{\lambda_1(p-\gamma_1)} {p\gamma_1}\|h_1\|_{L^{\frac{p}{p-\gamma_1}}(V)}^{\frac{p}{p-\gamma_1}}
    +  \frac{\lambda_2(q-\gamma_2)}{q\gamma_2}\|h_2\|_{L^{\frac{q}{q-\gamma_2}}(V)}^{\frac{q}{q-\gamma_2}}
   <  \frac{\alpha+\beta - \max\{p,q\}}{\alpha+\beta}
        M_1^{\frac{\alpha+\beta}{\alpha+\beta-\max\{p,q\}}}
        \left(\frac{\max\{p,q\}}{(\alpha+\beta)M_2}\right)^{\frac{\max\{p,q\}}{\alpha+\beta-\max\{p,q\}}}.
\end{cases}
\end{eqnarray}
\vskip2mm
\noindent
{\bf Theorem 5.1.} {\it Let $G=(V,E)$ be a finite graph. If $(\lambda_1,\lambda_2)$ satisfies (\ref{eqqq0}), then system (\ref{a2}) has at least one nontrivial solution of positive energy and one nontrivial solution of negative energy.

 }

\vskip2mm
\noindent
{\bf Theorem 5.2.} {\it Let $G=(V,E)$ be a finite graph.  For each $\lambda_1>0$, if $(u,v)$ is a semi-trivial solution of  system (\ref{a2}) and $(u,v)=(u,0)$, then
 $$
 \|u\|_{W^{m_1,p}(V)} \le \left(\lambda_1 H_1(V)C^{\gamma_1}_p(V)\right)^{\frac{1}{p-\gamma_1}},
 $$
 where $H_1(V)=\max_{x\in V} h_1(x)$.
 Similarly,  for each $\lambda_2>0$, if $(u,v)$ is a semi-trivial solution of  system (\ref{a2}) and $(u,v)=(0,v)$, then
 $$
 \|v\|_{W^{m_2,q}(V)} \le \left(\lambda_2 H_2(V)C^{\gamma_2}_q(V)\right)^{\frac{1}{q-\gamma_2}},
 $$
  where $H_2(V)=\max_{x\in V} h_2(x)$.
 }

\vskip2mm
 \par
 Denote
 \begin{eqnarray*}
 \lambda_0(V)   =   \frac{p-\gamma}{H_0(V)} C_p^{-\gamma}(V)
\left(\left(C_0(V) C_p^{\alpha}(V)\right)^{p-\alpha}(\alpha-p)^{\alpha-p}(\alpha-\gamma)^{\gamma-\alpha}\right)^{\frac{1}{p-\gamma}},
\end{eqnarray*}
\begin{eqnarray*}
\lambda_{\star}(V)  =   \frac{\gamma(\alpha-p)}{p\alpha H_0 C_p^\gamma(V)} \left(C_0(V)C^{\alpha}_p (V)\right)^{\frac{p- \gamma}{p-\alpha}},\ \
\lambda_{\star\star}(V)  =  \min\{\lambda_0(V), \lambda_{\star}(V)\},
\end{eqnarray*}
where $H_0(V)=\max_{x \in V} h(x)$, $C_0(V)=\max_{x \in V}c(x)$.
\vskip2mm
\noindent
{\bf Theorem 5.3.} {\it Let $G=(V,E)$ be a finite graph.  If $\lambda \in (0, \lambda_0(V))$, then equation (\ref{a3}) has at least one nontrivial solution of positive energy and one nontrivial solution of negative energy. Furthermore, if $\lambda \in (0, \lambda_{\star\star}(V))$, the solution of negative energy is the ground state solution of (\ref{a3}).}

\vskip3mm
 \noindent
\noindent{\bf Acknowledgments}

\noindent
 The authors sincerely thank the reviewer for his/her valuable comments which help us correct the proofs of Lemma 2.2 and Lemma 3.1, and improve the writing of the manuscript.

\vskip3mm
 \noindent
\noindent{\bf Funding information}

\noindent
This project is supported by Yunnan Fundamental Research Projects (grant No: 202301AT070465) and  Xingdian Talent Support Program for Young Talents of Yunnan Province.

\vskip2mm
\renewcommand\refname{References}
{}
\end{document}